\documentclass[11pt]{article}

\usepackage[a4paper,margin=1in]{geometry}
\usepackage{amsmath,amssymb,mathtools,bm}
\usepackage{graphicx}
\usepackage{booktabs}
\usepackage{array}
\usepackage{microtype}
\usepackage{xcolor}
\usepackage[round,authoryear]{natbib}
\usepackage[colorlinks=true, citecolor=blue, linkcolor=red, urlcolor=teal]{hyperref}
\usepackage{authblk}
\usepackage[title]{appendix}

\newcommand{\ud}{\ensuremath{\mathrm{d}}}

\newcommand{\bmi}{\ensuremath{\boldsymbol}}

\newcommand{\bnabla}{\ensuremath{\boldsymbol\nabla}}
\newcommand{\bcdot}{\ensuremath{\boldsymbol\cdot}}

\newcommand{\p}{\ensuremath{\partial}}
\newcommand{\beq}{\begin{equation}}
	\newcommand{\eeq}{\end{equation}}
\newcommand{\bal}{\begin{align}}
	\newcommand{\eal}{\end{align}}

\renewcommand{\arraystretch}{1.1}  

\title{\bfseries Weakly non-linear analysis of phototaxis-induced bioconvection in an isotropic scattering algal suspension under normal and oblique collimated irradiation}

\author[1]{S. K. Rajput}
\author[2]{A. Rathi}
\affil[1]{Department of Mathematics, School of Computer Science and Artificial Intelligence, SR University, Warangal 506371, India.}
\affil[2]{Department of Mathematics, PDPM Indian Institute of Information Technology Design and Manufacturing, Jabalpur 482005, India.}

\date{}

\begin{document}
\maketitle

\begin{abstract}
We investigate the weakly nonlinear evolution of phototaxis-induced bioconvection in an absorbing and isotropically scattering algal suspension subjected to oblique and normal collimated irradiation. Two hydrodynamic configurations are
considered, a rigid lower boundary with either a stress-free or a rigid upper surface. A multiple-scales expansion is carried to third order, and the Fredholm solvability condition yields a Stuart--Landau equation for the amplitude of the critical mode. Oblique irradiation affects the onset threshold in opposite ways in the two configurations. It progressively lowers the critical Rayleigh number beneath a stress-free surface and progressively raises it beneath a rigid surface. Beneath a stress-free surface the onset remains stationary, but the pitchfork bifurcation changes from supercritical to subcritical and becomes strongly subcritical at large incidence angles, where the preferred pattern wavelength also increases markedly. Beneath a rigid surface the critical mode changes from stationary to oscillatory at
moderate obliquity, giving a supercritical Hopf bifurcation, and returns to stationary at large obliquity, where the pitchfork bifurcation becomes subcritical.
\end{abstract}

\section{Introduction}
Light plays a dual role in suspensions of photosynthetic microorganisms:
it provides the environmental cue that directs cell swimming and, at the
same time, it is redistributed by the microorganisms through absorption
and scattering. Consequently, the optical field cannot generally be
regarded as an externally prescribed quantity. A change in cell
concentration modifies the local radiation field, the modified radiation
changes the swimming bias, and the resulting cell redistribution feeds
back into the fluid motion. This closed radiative--biological--hydrodynamic
feedback is one of the distinctive features of phototactic
bioconvection.


The fluid motion arises because motile microorganisms are typically
slightly denser than the surrounding liquid. Directed swimming can
therefore create strongly non-uniform cell distributions and, when a
cell-rich region overlies a more dilute region, the associated density
stratification may become gravitationally unstable. The resulting
large-scale collective motion is known as \emph{bioconvection}
\citep{13pedley1992,14hill2005,15bees2020,47Rajput2026BioconvectiveLS}. For phototactic
microorganisms, this instability is particularly sensitive to the
radiation field because the preferred swimming direction itself depends
on the local light intensity.

In a phototactic suspension, microorganisms may swim toward weakly
illuminated regions and reverse their swimming direction when the local
intensity exceeds a species-dependent critical value. The competition
between phototactic swimming and cell diffusion can then produce an
interior concentrated sublayer. The region below this sublayer may be
gravitationally unstable, while the region above it is comparatively
stable, giving rise to penetrative phototactic bioconvection
\citep{8vincent1996,6ghorai2005}. When scattering is present, this
mechanism becomes more involved because the cells influence not only
the attenuation of the incident beam but also the redistribution of
radiation over different propagation directions.

The classical continuum theory of phototactic bioconvection was first
developed for a vertically illuminated, absorbing and non-scattering
suspension \citep{8vincent1996}. Subsequent studies incorporated
penetrative convection and radiative scattering
\citep{6ghorai2005,4ghorai2010,23ghorai2013}. In an isotropically
scattering suspension, the radiation field contains a directly
transmitted collimated component and a diffuse component generated by
scattering. The total light intensity and radiative flux are therefore
determined by both the attenuation of the incident beam and its
redistribution by the microorganisms. This coupling can substantially
alter the location and structure of the basic cell-concentration
sublayer and hence the stability characteristics of the suspension.

A further control parameter is the direction of illumination. In many
natural and laboratory situations the incident beam is not perpendicular
to the upper surface. We therefore describe the irradiation by an
external incidence angle $\theta_i$, with the normal-incidence
configuration recovered continuously at $\theta_i=0^\circ$. Finite
values of $\theta_i$ correspond to oblique collimated irradiation and
modify the path of the direct beam through the suspension. In an
isotropically scattering medium, this change in the collimated
contribution also modifies the diffuse radiation field generated by
scattering. Linear studies of related phototactic systems have shown
that varying the incidence angle can change the basic concentration
profile, the critical Rayleigh number and wavenumber, and even the
stationary or oscillatory character of the most unstable mode
\citep{32rajput2025impact}.

These linear predictions, however, describe only the loss of stability
of the basic state. They do not determine the fate of the disturbance
once its amplitude is no longer infinitesimal. In particular, the
critical Rayleigh number alone cannot establish whether the emerging
branch is locally supercritical or subcritical, whether nonlinear
interactions saturate the disturbance, or how the frequency of an
oscillatory instability changes after onset. In the present scattering
problem these questions are especially significant because the
perturbed concentration modifies both the collimated and diffuse
radiation fields, thereby introducing an additional nonlinear feedback
into the phototactic cell flux.

To resolve this near-onset dynamics, we employ a weakly nonlinear
multiple-scales analysis. The Rayleigh number is expanded about its
critical value and the amplitude of the critical mode is allowed to
evolve on a slow time scale. The first-order problem recovers the
critical linear eigenmode. Quadratic interactions at second order
generate a correction to the mean state and higher spatial harmonics,
while the third-order resonant problem is constrained by a Fredholm
solvability condition. This leads to a Stuart--Landau amplitude equation
\citep{44Malkus1958FiniteAC,45Stuart1960OnTN,46Watson1960OnTN}, whose
cubic coefficient determines the local nonlinear saturation and
bifurcation character. For an oscillatory critical mode, the coefficient
is complex and also gives the leading nonlinear phase and frequency
correction.

The present study applies this framework to an absorbing and
isotropically scattering algal suspension subjected to oblique
collimated irradiation. The incidence angle is treated as a continuous
parameter, so that normal illumination appears naturally as the
$\theta_i=0^\circ$ limit of the same model. Two mechanical
configurations are examined: a rigid lower boundary with either a
stress-free or a rigid upper surface. The aim is to determine how the
irradiation geometry and isotropic scattering
enter the weakly nonlinear evolution of stationary and oscillatory
bioconvective modes near the instability threshold.
\section{Mathematical formulation}\label{formulation}
We consider an infinitely extended horizontal layer of fluid of depth $H$, populated by phototactic microorganisms, in a Cartesian coordinate system $Oxyz$ with the $z$-axis directed vertically upwards. The suspension is illuminated from above by a collimated light beam incident at a prescribed oblique angle $\theta_i$ to the vertical. The plane of incidence is chosen to coincide with the $xz$-plane. On entering the suspension, the beam is refracted through an angle \(\theta_0\), determined by Snell's law,
$\sin\theta_i = 1.333\,\sin\theta_0$, where \(1.333\) is the refractive index of water relative to air. Let \(I(\bmi{x},\bmi{s})\) denote the radiative intensity at position \(\bmi{x}=(x,y,z)\) in the propagation direction $\bmi{s}$. In spherical polar coordinates,
$\bmi{s}=\sin\theta\cos\phi\,\hat{\bmi{i}}+\sin\theta\sin\phi\,\hat{\bmi{j}}+\cos\theta\,\hat{\bmi{k}}$, where \(\theta\) is the polar angle measured from the positive $z$-axis and $\phi$ is the azimuthal angle measured from the positive $x$-axis in the $xy$-plane.

\begin{figure}[!htbp]
	\centering
	\includegraphics[scale=0.5]{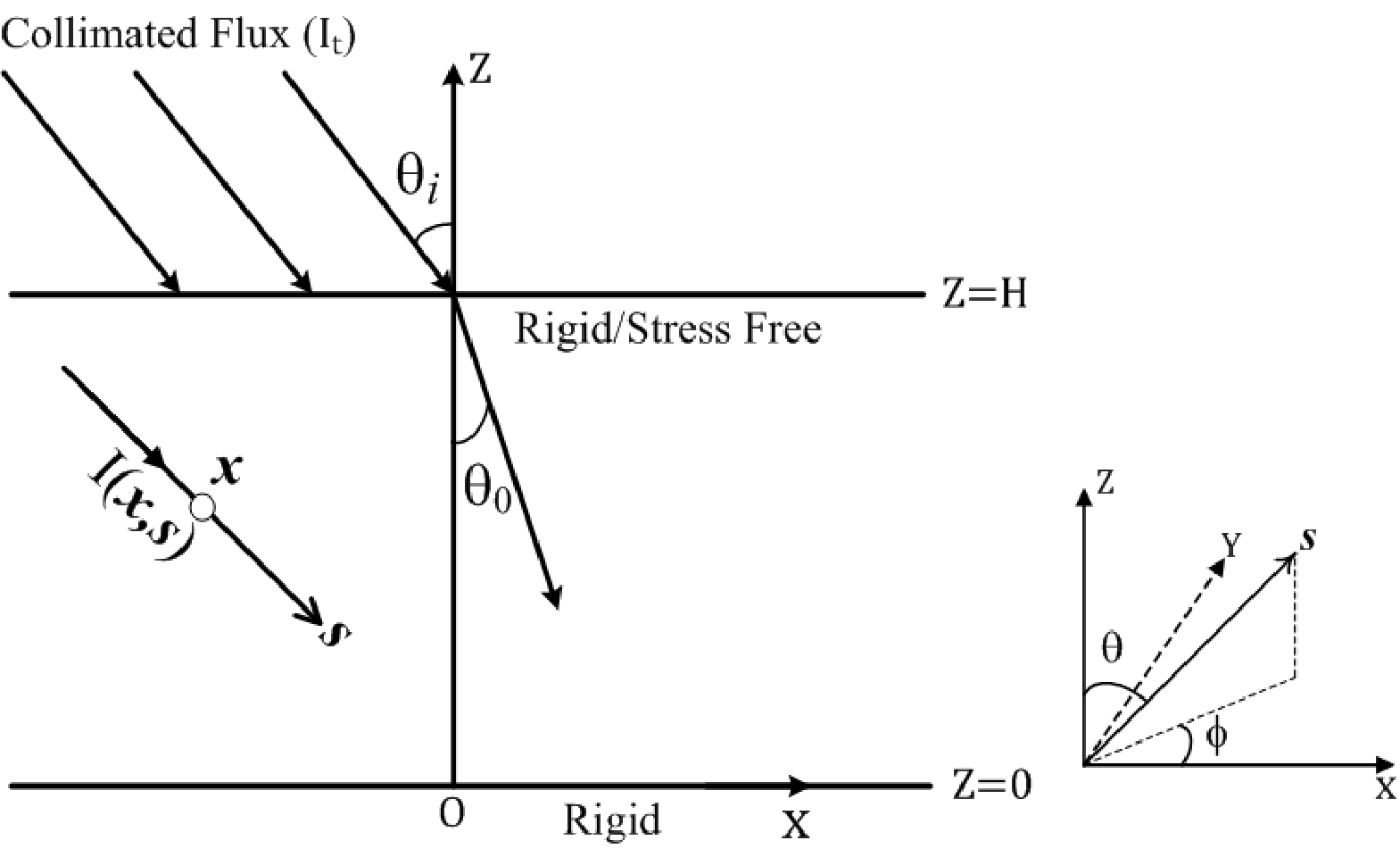}
	\caption{\footnotesize{Obliquely incident collimated irradiation on the the upper surface of an absorbing and isotropically scattering suspension.}}
	\label{fig2}
\end{figure}

\subsection{Governing equations}\label{GE}
Following classical continuum models of bioconvection~\citep{13pedley1992}, the cells are treated as a monodisperse population described by the concentration field \(n(\boldsymbol{x},t)\), while \(\boldsymbol{u}\) denotes the mean fluid velocity. The suspension is assumed dilute and incompressible, and the Boussinesq approximation is adopted. The governing equations are therefore
\begin{equation}
	\bmi{\nabla} \bcdot \bmi{u} = 0.
\end{equation}
\begin{align}
	\rho\left(\frac{\partial\bmi{u}}{\partial t}+\bmi{u}\bcdot\bmi{\nabla}\bmi{u}\right)= -\bmi{\nabla} P_e + \mu \bmi{\nabla}^2 \bmi{u} - n v g \delta \rho \hat{\bmi{k}}.
\end{align}
where $\rho$ and $\mu$ are the density and dynamic viscosity of the fluid, respectively, \(v\) is the volume of an individual cell and $\delta\rho$ is the excess cell density relative to the surrounding fluid. \\
Conservation of cells satisfies
\begin{equation}
	\frac{\partial n}{\partial t} + \bmi{\nabla} \cdot \bmi{F} = 0,
\end{equation}
where the cell flux is
$$
\bmi{F}=n \bmi{u} + n \bmi{U_c} - \bmi{D} \bcdot \bmi{\nabla} n.
$$
The cell-flux formulation assumes purely phototactic swimming, with viscous torques associated with horizontal motion neglected. Following \cite{6ghorai2005}, cell diffusion is taken to be isotropic and constant, such that
$\bmi{D}=D\bmi{I}$, where $D$ denotes the cell diffusivity.


\subsection{Average swimming orientation}\label{mean direction}
The radiative intensity \(I(\bmi{x},\bmi{s})\) within an absorbing and scattering medium is governed by the radiative transfer equation (RTE)~\citep{33modest2003}

\begin{equation}
	\frac{\ud I(\bmi{x},\bmi{s})}{\ud s} + (\kappa_a + \kappa_s)I(\bmi{x},\bmi{s}) = \frac{\kappa_s}{4\pi}\int_{0}^{4\pi}I(\bmi{x},\bmi{s}') f(\bmi{s}',\bmi{s}) d\Omega'
\end{equation} 
where \(\kappa_a\) and \(\kappa_s\) are the absorption and scattering coefficients, respectively, and \(\Omega'\) is the solid angle. Here, \(f(\bmi{s}',\bmi{s})\) is the scattering phase function, describing the probability that a radiation traveling in the direction \(\bmi{s}'\) is scattered into the direction \(\bmi{s}\). For simplicity, the phase function \(f(\bmi{s}',\bmi{s}) = 1\) for isotropic scattering~\citep{33modest2003}.
At the top boundary, the intensity of light at a spatial location \(\bmi{x_b} = (x,y,z=H)\) is 
\begin{equation}
	I(\bmi{x_b},\bmi{s}) = \mathrm{I_t} \delta (\bmi{s}-\bmi{s}_0)
\end{equation}
where \(\mathrm{I_t}\) is the magnitude of the oblique incident collimated sunlight. Also \(\bmi{s}_0 = \sin\theta_0 \cos \phi_0 \,\hat{\bmi{i}} + \sin\theta_0 \sin \phi_0 \,\hat{\bmi{j}} + \cos\theta_0 \,\hat{\bmi{k}}\) is the incident direction. Moreover, \(\delta (\bmi{s}-\bmi{s}_0)\) is the Dirac delta function, which satisfies at \(\bmi{s}=\bmi{s}_0\)
\begin{equation}
	\int_{0}^{4\pi} g(\bmi{s}) \delta (\bmi{s}-\bmi{s}_0) d \Omega = g(\bmi{s}_0)
\end{equation}

Here, the absorption and scattering are considered as linearly dependent on the cell concentration so that \(\kappa_a = \varkappa n(\bmi{x})\) and \(\kappa_s = \varsigma n(\bmi{x})\). Using the single scattering albedo \(\omega = \kappa_s/(\kappa_a + \kappa_s)\), the RTE can be written as
\begin{equation}
	\frac{\ud I(\bmi{x},\bmi{s})}{\ud s} + \alpha I(\bmi{x},\bmi{s}) = \frac{\omega \alpha }{4\pi}\int_{0}^{4\pi}I(\bmi{x},\bmi{s}') d\Omega'
\end{equation} 
where, \(\alpha = (\varkappa+\varsigma)n\) is the extinction coefficient. In addition, the parameter \(\omega \in [0,1]\), where \(\omega=0\) signifies a completely absorbing suspension, and \(\omega=1\) corresponds to a purely scattering suspension~\citep{23ghorai2013}. The intensity of light \(G(\bmi{x})\) and the radiative heat flux \(\bmi{q(x)}\) are defined as follows~\citep{25panda2020}
\begin{equation}
	G(\bmi{x}) = \int_{0}^{4\pi} I(\bmi{x},\bmi{s}) d \Omega \quad \text{and} \quad \bmi{q(x)} = \int_{0}^{4\pi} I(\bmi{x},\bmi{s}) \bmi{s} d \Omega.
\end{equation} 

For many microbial species, their swimming speed doesn't depend on the illumination, position, time, and direction~\citep{10hill1997}. Let \(\bmi{p}\) and \( U_c \) be the mean swimming direction and mean swimming speed of cells. Thus, the mean swimming velocity \(\bmi{U_c}\) is 
\begin{equation}
	\bmi{U_c} = U_c \bmi{p}.
\end{equation}

Here, the mean swimming direction, \(\bmi{p}\), is~\citep{25panda2020}
\begin{equation}
	\bmi{p} = -T(G) \frac{\bmi{q}}{|\bmi{q}|},
\end{equation}
where \( T(G) \) is the phototaxis function and given by
\begin{equation}
	T(I) =
	\begin{cases} 
		\geq 0, & G(\bmi{x}) \leq G_c, \\
		< 0, & G(\bmi{x}) > G_c.
	\end{cases}
\end{equation}

Here, \(T(G)\) represents a light-intensity-dependent function whose form varies with the species of microorganisms under consideration~\citep{8vincent1996}. Here, the phototaxis function \(T(G)\) is modelled as~\citep{32rajput2025impact}
\begin{equation}\label{13}
	T(G) = 0.8 \sin \left[\frac{3\pi}{2} \chi(G) \right] - 0.1 \sin \left[\frac{\pi}{2} \chi(G) \right],
\end{equation}
where \( \chi(G) = \frac{G}{2.5} \, e^{\{0.32\,(2.5-G)\}} \) with critical light intensity \(G_c = 1\).

\subsection{Boundary conditions}
At the lower boundary $z = 0$, rigid no-slip condition is imposed while for upper boundary $z=H$, rigid and stress-free both case are considered here. In addition, the algal cells can not flow through horizontal boundaries. Accordingly, the boundary conditions for the governing system are 
\begin{equation}
	\bmi{u} \cdot \hat{\bmi{k}} = \bmi{F} \cdot \hat{\bmi{k}} = 0 \quad \text{at} \quad z = 0,\, H
\end{equation}
For a rigid boundary
\begin{equation}
	\bmi{u} \times \hat{\bmi{k}} = 0 \quad \text{at} \quad z = 0,\, H
\end{equation}		 	
whilst for a stress-free boundary
\begin{equation}
	\frac{\p^2}{\p z^2} (\bmi{u} \cdot \hat{\bmi{k}}) = 0 \quad \text{at} \quad z = H
\end{equation}
The upper surface is illuminated by a uniform obliquely incident collimated solar flux of constant magnitude \(\mathrm{I_t}\). Since both the upper and lower boundaries are assumed to be non-reflecting, the corresponding boundary conditions are
\begin{equation}
	\left.
	\begin{aligned}
		I(x,y,z=H, \theta, \phi) &= \mathrm{I_t} \delta (\bmi{s}-\bmi{s}_0), && \pi/2 \leq \theta \leq \pi, \\
		I(x,y,z=0, \theta, \phi) &= 0, && 0 \leq \theta \leq \pi/2
	\end{aligned}
	\right\}
\end{equation}

\subsection{Scaling of the Equations}

Here, length is scaled by the suspension depth \(H\), time by the diffusive time scale \(H^2/D\), and the bulk fluid velocity by \(D/H\). The cell concentration is scaled by \(N_0\), and the pressure by \(\mu D/H^2\)~\citep{6ghorai2005,25panda2020}. In terms of the scaled variables, the governing equations are 
\begin{equation}
	\bmi{\nabla} \bcdot \bmi{u}=0
\end{equation}	
\begin{equation}
	S_c^{-1}\left(\frac{\partial{\bmi{u}}}{\partial{t}}+\left(\bmi{u}\bcdot\bmi{\nabla}\right)\bmi{u}\right)=-\bmi{\nabla} P_e+{\bmi{\nabla}}^2\bmi{u}-R_a \, n\, \hat{\bmi{k}},
\end{equation}	
\begin{equation}
	\frac{\partial n}{\partial t} + \boldsymbol{\nabla} \bcdot \boldsymbol{F} = 0,
\end{equation}
where
\begin{equation*}
	\boldsymbol{F} = n \boldsymbol{u} +  n V_c \boldsymbol{p} -  \boldsymbol{\nabla} n.
\end{equation*}	

Here \( S_c = \frac{\mu}{\rho D} \) is the Schmidt number, \( V_c = \frac{U_c H}{D} \) is scaled swimming speed and \(R_a = \frac{N_0 \vartheta g \Delta \rho H^{3}}{\mu D}\) is the Rayleigh number. In terms of the non-dimensional variables, RTE takes the form
\begin{equation}
	\frac{\ud I(\bmi{x},\bmi{s})}{\ud s} + \alpha_H I(\bmi{x},\bmi{s}) = \frac{\omega \alpha_H n }{4\pi}\int_{0}^{4\pi}I(\bmi{x},\bmi{s}') d\Omega'
\end{equation}
where \(\alpha_H = (\varkappa + \varsigma) N_0 H \) is the scaled extinction coefficient. Using the direction cosines \((\xi_x,\xi_y,\xi_z) = (\sin\theta \cos\phi, \sin\theta \sin\phi, \cos\theta)\) of the unit vector \(\boldsymbol{s}\), the RTE becomes~\citep{33modest2003} 
\begin{align}
	\xi_x \frac{\partial I}{\partial x} + \xi_y \frac{\partial I}{\partial y} + \xi_z \frac{\partial I}{\partial z} + \alpha_H n I(\bmi{x},\bmi{s}) = \frac{\omega \alpha_H n }{4\pi}\int_{0}^{4\pi}I(\bmi{x},\bmi{s}') d\Omega'
\end{align}	
After scaling, the boundary conditions are 	
\begin{equation}
	\bmi{u} \cdot \hat{\bmi{k}} = \bmi{F} \cdot \hat{\bmi{k}} = 0 \quad \text{at} \quad z = 0,\, 1
\end{equation}
For a rigid boundary
\begin{equation}
	\bmi{u} \times \hat{\bmi{k}} = 0 \quad \text{at} \quad z = 0,\, 1
\end{equation}		 	
whilst for a stress-free boundary
\begin{equation}
	\frac{\p^2}{\p z^2} (\bmi{u} \cdot \hat{\bmi{k}}) = 0 \quad \text{at} \quad z = 1
\end{equation}
Due to the scaling of light intensity by \(\mathrm{I_t}\), the corresponding boundary conditions are 
\begin{equation}
	\left.
	\begin{aligned}
		I(x,y,z=1, \theta, \phi) &= \delta (\bmi{s}-\bmi{s}_0), && \pi/2 \leq \theta \leq \pi, \\
		I(x,y,z=0, \theta, \phi) &= 0, && 0 \leq \theta \leq \pi/2
	\end{aligned}
	\right\}
\end{equation}

\section{The basic solution}
Eqs.(17)-(19) and Eq.(21) subject to the boundary conditions have a steady state solution such that
\begin{equation}
	\begin{pmatrix}
		\bmi{u} \\ n \\ I
	\end{pmatrix}
	=
	\begin{pmatrix}
		0  \\ n_0(z)  \\ I_0(z,\theta)  
	\end{pmatrix},
\end{equation}	

At the basic state, the total light intensity \(G_0(z)\) and the radiative heat flux \(\bmi{q}_0(z)\) are given by
\begin{equation}
	G_0(z) = \int_{0}^{4\pi} I_0(z,\theta) d \Omega \quad \text{and} \quad \bmi{q}_0(z) = \int_{0}^{4\pi} I_0(z,\theta) \bmi{s} d \Omega.
\end{equation}

Since \(I_0(z,\theta)\) is independent of the azimuthal angle \(\phi\), so that the horizontal component of \(\bmi{q}_0\) is zero. Therefore, \(\bmi{q}_0 = - q_0 \hat{\bmi{k}}\), where \(q_0 = |\bmi{q}_0|\).

The basic state radiative intensity, \(I_0(z,\theta)\), is governed by
\begin{equation}
	\xi_z \frac{\partial I_0}{\partial z} + \alpha_H n_0 I_0(z,\theta) = \frac{\omega \alpha_H n_0}{4\pi} G_0(z) 
\end{equation}
The basic state radiative intensity \(I_0(z,\theta)\) can be expressed as: \(I_0(z,\theta) = I_0^c(z,\theta)+I_0^d(z,\theta)\). The collimated part $I_0^c(z,\theta)$ satisfies 
\begin{align}
	\frac{\partial I_0^{c}(z,\theta)}{\partial z} + \frac{\alpha_H n_0(z)}{\xi_z} I_0^{c}(z,\theta) = 0
\end{align}
with boundary condition
\begin{equation}
	I_0^{c}(z=1,\theta) = \delta (\bmi{s}-\bmi{s}_0), \quad \pi/2 \leq \theta \leq \pi
\end{equation}
On solving Eqs.(29) and (30), we get
\begin{equation}
	I_0^{c} = \exp\left(\int_{z}^{1}\frac{\alpha_H n_0(z')}{\cos\theta} dz'\right) \delta (\bmi{s}-\bmi{s}_0)
\end{equation}
The basic state diffuse intensity \(I_0^{d}(z,\theta)\) satisfies 
\begin{equation}
	\frac{\partial I_0^d}{\partial z} + \frac{\alpha_H n_0}{\xi_z} I_0^d(z,\theta) = \frac{\omega \alpha_H n_0}{4\pi \xi_z} G_0(z)
\end{equation}
with boundary conditions 
\begin{equation}
	\left.
	\begin{aligned}
		I_0^d(z=1,\theta) &= 0, && \pi/2 \leq \theta \leq \pi,\  \\
		I_0^d(z=0,\theta) &= 0, && 0 \leq \theta \leq \pi/2,\ 
	\end{aligned}
	\right\}
\end{equation}
The total intensity \(G_0(z)\) is expressed as the sum of diffuse and collimated parts 
\begin{equation}
	G_0(z) = G_0^c(z) + G_0^d(z)
\end{equation}
where
\begin{equation*}
	G_0^{c} = \int_{0}^{4\pi} I_0^{c}(z,\theta) d\Omega = \exp\left(-\frac{\alpha_H}{\cos \theta_0}\int_{z}^{1} n_0(z') dz'\right)
\end{equation*}
and
\begin{equation*}
	G_0^{d} = \int_{0}^{4\pi} I_0^{d}(z,\theta) d\Omega
\end{equation*}
Introducing a new variable, \(\tau\), such that 
\begin{equation}
	\tau = \int_{z}^{1}\alpha_H n_0(z') dz'
\end{equation}
therefore, the total intensity $G_0$ is given by
\begin{equation}
	\begin{aligned}
		G_0(\tau) = e^{-\tau/\cos \theta_0} +  \frac{\omega}{2}\int_{0}^{\alpha_H} G_0(\tau')E_1(|\tau-\tau'|) d\tau'
	\end{aligned}
\end{equation}		
where, \(E_1(x)\) is the exponential integral function. Eq.(36) is solved using the subtraction of singularity technique~\citep{19crosbie1985application,35press1992}. \\
The radiative heat flux is given by
\begin{equation}
	\bmi{q}_0(z) = \bmi{q}_0^c(z) + \bmi{q}_0^d(z)
\end{equation}
where
\begin{equation*}
	\bmi{q}_0^{c} = \int_{0}^{4\pi} I_0^{c}(z,\theta) \bmi{s} d\Omega = -\cos \theta_0\exp\left(-\frac{\alpha_H}{\cos \theta_0}\int_{z}^{1} n_0(z') dz'\right)\hat{\bmi{k}}
\end{equation*}
and
\begin{equation*}
	\bmi{q}_0^{d} = \int_{0}^{4\pi} I_0^{d}(z,\theta) \bmi{s} d\Omega
\end{equation*}
Therefore, the mean swimming direction of cells in the basic state becomes
\begin{equation}
	\bmi{p}_0 = -T(G_0) \frac{\bmi{q}_0}{q_0} = T(G_0)\hat{\bmi{k}}.
\end{equation}
The basic state concentration \( n_0(z) \) satisfies the equation	
\begin{equation}
	\frac{\ud n_0}{\ud z} = V_c \,T_0 \,n_0, \quad \text{where} \quad T_0 = T(G_0)
\end{equation}	
Here, Eq.(39) is supplemented with the conservation relation 
\begin{equation}
	\int_0^1 n_0(z) dz = 1
\end{equation}	
Next, Eqs.(39)-(40) constitute a boundary value problem (BVP) which is numerically solved using a shooting method~\citep{35press1992}. 


\section{Perturbations}

The basic state solution is perturbed as 	
\begin{align*}
	&\bmi{u} = 0 + \epsilon \bmi{u}_1 + \epsilon^2 \bmi{u}_2 + \epsilon^3 \bmi{u}_3 + \cdots\\
	&n = n_0 + \epsilon n_1 + \epsilon^2 n_2 + \epsilon^3 n_3 + \cdots\\
	&I = I_0 + \epsilon I_1 + \epsilon^2 I_2 + \epsilon^3 I_3 + \cdots\\
	&\bmi{p} = \bmi{p}_0 + \epsilon \bmi{p}_1 + \epsilon^2 \bmi{p}_2 + \epsilon^3 \bmi{p}_3 + \cdots\\
	&R_a = R_c + \epsilon^2 R_2 + \epsilon^4 R_4 + \cdots\\
	&P_e = P_{e_0} + \epsilon P_{e_1} + \epsilon^2 P_{e_2} + \epsilon^3 P_{e_3} + \cdots
\end{align*}	 	

The fast and slow times are
$\tau=t$ and $s=\epsilon^2t$, so
\begin{equation*}
	\frac{\partial}{\partial t} = \frac{\partial}{\partial \tau} + \epsilon^2 \frac{\partial}{\partial s}.
\end{equation*}


\subsection{First Order System}

At $O(\epsilon)$, the governing equations are
\begin{equation}
	\boldsymbol{\nabla} \bcdot \boldsymbol{u}_1 = 0,
\end{equation}
\begin{equation}
	S_c^{-1} \frac{\partial \bmi{u}_1}{\partial \tau} = -\boldsymbol{\nabla} P_{e_1} + \boldsymbol{\nabla}^{2} \boldsymbol{u}_1 - R_c \, n_1 \, \hat{\boldsymbol{k}},
\end{equation}	
\begin{equation}
	\frac{\partial n_1}{\partial \tau} + \bmi{\nabla} \bcdot [n_0 \bmi{u}_1+ V_c (n_0 \bmi{p}_1 + n_1 \bmi{p}_0)] - \bmi{\nabla}^2 n_1 = 0
\end{equation}
Now, we are taking the curl twice of Eq.(42) to remove the pressure term and retaining the vertical component of $\bmi{u}_1$ by removing the horizontal components of the result. Thus, the governing equations become
\begin{equation}
	\bnabla^4 w_1 - S_c^{-1} \frac{\partial}{\partial \tau} \Bigl(\bnabla^2 w_1\Bigr) - R_c \left(\frac{\partial^2}{\partial x^2} + \frac{\partial^2}{\partial y^2}\right) n_1 = 0
\end{equation}
\begin{equation}
	\frac{\partial n_1}{\partial \tau} + w_1 \frac{\p n_0}{\p z} + V_c \frac{\p}{\p z} \Big(T_0 \, n_1 + n_0 T'(G_0) G_1 \Big) - \frac{V_c n_0 T_0}{q_0} \left(\frac{\p q_1^{x}}{\p x}+\frac{\p q_1^{y}}{\p y}\right)  = \bmi{\nabla}^2 n_1
\end{equation}
We define a new variable $\widehat{n}_1=\int_{z}^{1} n_1 \ud z$, the governing system can be written as follows
\begin{equation}
	\bmi{\mathcal{L}} \bmi{X}_1 = \bmi{F}_1
\end{equation}
\begin{equation}
	\implies
	{\setlength{\arraycolsep}{12pt}
		\renewcommand{\arraystretch}{1.6}
		\begin{bmatrix}
			\bnabla^4 - S_c^{-1} \dfrac{\partial}{\partial \tau} \bnabla^2 & R_c \left(\dfrac{\partial^2}{\partial x^2} + \dfrac{\partial^2}{\partial y^2}\right) \dfrac{\partial}{\partial z} \\
			-\dfrac{\partial n_0}{\partial z} & \dfrac{\partial^2}{\partial \tau \partial z} - \mathcal{F}
	\end{bmatrix}}
	\begin{bmatrix}
		w_1 \\
		\widehat{n}_1
	\end{bmatrix}
	=
	\begin{bmatrix}
		f_{11} \\
		f_{12}
	\end{bmatrix}.
\end{equation}
where 
$$
f_{11} = 0
$$
$$
f_{12} = V_c \frac{\partial}{\partial z}\left(n_0 T'(G_0) G_1^d \right)-\frac{V_c n_0 T_0}{q_0} \left(\frac{\p q_1^{x}}{\p x}+\frac{\p q_1^{y}}{\p y}\right)
$$
and $$ \mathcal{F} = \bnabla^2 \frac{\p}{\p z} - V_c \left[T_0 \frac{\p^2}{\p z^2} + T'(G_0)\left(2 (\alpha_H/\cos\theta_0) n_0 G_0^c + \frac{\p G_0^d}{\p z}\right)\frac{\p}{\p z} + (\alpha_H/\cos\theta_0) \frac{\partial}{\partial z}\Big(n_0 G_0^c T'(G_0)\Big)\right] $$
Applying normal mode analysis
$$
w_1=A(s)W_1(z)\exp\{i(lx+my) + \gamma \tau \} + c.c. ,\quad \widehat{n}_1=A(s)\widehat{\Theta}_1(z) \exp\{i(lx+my) + \gamma \tau \} + c.c.,
$$
The governing system, rewritten using \( D \equiv d/dz \) take the following form 
\begin{equation}
	{\setlength{\arraycolsep}{8pt}
		\renewcommand{\arraystretch}{1.2}
		\begin{bmatrix}
			(D^2-k^2)^2-\gamma S_c^{-1}(D^2-k^2) & -k^2 R_c D \\
			-D n_0 & \mathcal{F}_1
	\end{bmatrix}}
	\begin{bmatrix}
		W_1 \\
		\widehat{\Theta}_1
	\end{bmatrix}
	=
	\begin{bmatrix}
		\Gamma_{11} \\
		\Gamma_{12}
	\end{bmatrix}.
\end{equation}
Here
$$ \Gamma_{11} = 0, \quad \Gamma_{12} = V_c \frac{\partial}{\partial z}\left(n_0 T'(G_0) \mathcal{G}_1^d \right)-\frac{V_c n_0 T_0}{q_0} \Big(l \mathcal{P}+m \mathcal{Q}\Big)$$
and
$$ \mathcal{F}_1 = 	D^3 - V_c T_0 \, D^2 -\left\{\gamma + k^2 + \Gamma_1(z)\right\} D - \Gamma_2(z) $$
where
\begin{align*} 
	&\Gamma_1(z) = V_c T'(G_0)\Big\{2 (\alpha_H/\cos\theta_0) n_0 G_0^c + D G_0^d\Big\}, \\
	&\Gamma_2(z) = V_c(\alpha_H/\cos\theta_0) D\Big(n_0 G_0 T'(G_0)\Big)
\end{align*}
with rigid boundary conditions
\begin{equation}
	W_1 = DW_1 = D^2 \widehat{\Theta}_1 - V_c T_0 \, D\widehat{\Theta}_1 - V_c T'(G_0) \mathcal{G}_1 \widehat{\Theta}_1 = 0 \quad \text{at} \quad z = 0, 1
\end{equation}	
while for a free upper boundary
\begin{equation}
	W_1 = D^2 W_1 = D^2 \widehat{\Theta}_1 - V_c T_0 \, D\widehat{\Theta}_1 = 0 \quad \text{at} \quad z = 1.
\end{equation}	
and
\begin{equation}
	\widehat{\Theta}_1 = 0 \quad \text{at} \quad z = 1.
\end{equation}
Here, $k=\sqrt{l^{2}+m^{2}}$ denotes the dimensionless horizontal wavenumber. Eqn~(48), together with the associated boundary conditions, defines an eigenvalue problem for the complex growth rate $\gamma$.

\subsection{Second Order System}
At $O(\epsilon^2)$, the governing system is as follows	
\begin{equation}
	\bnabla^4 w_2 - S_c^{-1} \frac{\partial}{\partial \tau} \biggl(\bnabla^2 w_2\biggr) - R_c \left(\frac{\partial^2}{\partial x^2} + \frac{\partial^2}{\partial y^2}\right) n_2 = -S_c^{-1}\Bigl[\bnabla \{\bnabla \bcdot (\bmi{u}_1 \bcdot \bnabla \bmi{u}_1)\} - \bnabla^2 (\bmi{u}_1 \bcdot \bnabla \bmi{u}_1) \Bigr]\bcdot \hat{\bmi{k}}
\end{equation}
\begin{equation}
	\frac{\partial n_2}{\partial \tau} + \bmi{\nabla} \bcdot (n_0 \bmi{u}_2 + n_1 \bmi{u}_1) + V_c \bmi{\nabla} \bcdot (n_0 \bmi{p}_2 + n_1 \bmi{p}_1 + n_2 \bmi{p}_0) - \bmi{\nabla}^2 n_2 = 0
\end{equation}
with rigid boundary conditions
\begin{subequations}
	\begin{equation}
		w_2 = \frac{\partial w_2}{\partial z} = \Big[ \bnabla n_2 - V_c \left(n_0 \bmi{p}_2 + n_1 \bmi{p}_1 + n_2 \bmi{p}_0 \right) \Big] \bcdot \hat{\bmi{k}} = 0, \quad \text{at} \quad z = 0, 1
	\end{equation}
	while for a free upper boundary
	\begin{equation}
		w_2 = \frac{\partial^2 w_2}{\partial z^2} = \Big[ \bnabla n_2 - V_c \left(n_0 \bmi{p}_2 + n_1 \bmi{p}_1 + n_2 \bmi{p}_0 \right) \Big] \bcdot \hat{\bmi{k}} = 0, \quad \text{at} \quad z = 1
	\end{equation}
\end{subequations}
Using a new variable $\widehat{n}_2=\int_{z}^{1} n_2 \ud z$, the governing system can be written as
\begin{equation}
	\bmi{\mathcal{L}} \bmi{X}_2 = \bmi{F}_2
\end{equation}
\begin{equation}
	\implies
	{\setlength{\arraycolsep}{12pt}
		\renewcommand{\arraystretch}{1.6}
		\begin{bmatrix}
			\bnabla^4 - S_c^{-1} \dfrac{\partial}{\partial \tau} \bnabla^2 & R_c \left(\dfrac{\partial^2}{\partial x^2} + \dfrac{\partial^2}{\partial y^2}\right) \dfrac{\partial}{\partial z} \\
			-\dfrac{\partial n_0}{\partial z} & \dfrac{\partial^2}{\partial \tau \partial z} - \mathcal{F}
	\end{bmatrix}}
	\begin{bmatrix}
		w_2 \\
		\widehat{n}_2
	\end{bmatrix}
	=
	\begin{bmatrix}
		f_{21} \\
		f_{22}
	\end{bmatrix}.
\end{equation}
where 
$$
f_{21} = -S_c^{-1}\Bigl[\bnabla \{\bnabla \bcdot (\bmi{u}_1 \bcdot \bnabla \bmi{u}_1)\} - \bnabla^2 (\bmi{u}_1 \bcdot \bnabla \bmi{u}_1) \Bigr]\bcdot \hat{\bmi{k}}
$$
$$
f_{22} = \bmi{\nabla} \bcdot (n_1 \bmi{u}_1) + V_c \bmi{\nabla} \bcdot (n_1 \bmi{p}_1)
$$

\subsection{Third Order System}
At $O(\epsilon^3)$, the governing system is as follows
\begin{align}
	\nonumber	\bnabla^4 w_3 - S_c^{-1} \frac{\partial}{\partial \tau} \biggl(\bnabla^2 w_3\biggr) - R_c \left(\frac{\partial^2}{\partial x^2} + \frac{\partial^2}{\partial y^2}\right) n_3 = S_c^{-1} \frac{\partial}{\partial s} \biggl(\bnabla^2 w_1\biggr) + R_2 \left(\frac{\partial^2}{\partial x^2} + \frac{\partial^2}{\partial y^2}\right) n_1 \\ - S_c^{-1}\Bigl[\bnabla \{\bnabla \bcdot (\bmi{u}_1 \bcdot \bnabla \bmi{u}_2)\} - \bnabla^2 (\bmi{u}_1 \bcdot \bnabla \bmi{u}_2) + \bnabla \{\bnabla \bcdot (\bmi{u}_2 \bcdot \bnabla \bmi{u}_1)\} - \bnabla^2 (\bmi{u}_2 \bcdot \bnabla \bmi{u}_1) \Bigr]\bcdot \hat{\bmi{k}}
\end{align}	
\begin{equation}
	\frac{\partial n_3}{\partial \tau} + \frac{\partial n_1}{\partial s} + \bmi{\nabla} \bcdot (n_0 \bmi{u}_3 + n_1 \bmi{u}_2 + n_2 \bmi{u}_1) + V_c \bmi{\nabla} \bcdot (n_0 \bmi{p}_3 + n_1 \bmi{p}_2 + n_2 \bmi{p}_1 + n_3 \bmi{p}_0) - \bmi{\nabla}^2 n_3 = 0
\end{equation}
with rigid boundary conditions
\begin{subequations}
	\begin{equation}
		w_3 = \frac{\partial w_3}{\partial z} = \Big[ \bnabla n_3 - V_c (n_0 \bmi{p}_3 + n_1 \bmi{p}_2 + n_2 \bmi{p}_1 + n_3 \bmi{p}_0)\Big]\bcdot \hat{\bmi{k}} = 0, \quad \text{at} \quad z = 0, 1
	\end{equation}
	while for a free upper boundary
	\begin{equation}
		w_3 = \frac{\partial^2 w_3}{\partial z^2} = \Big[ \bnabla n_3 - V_c (n_0 \bmi{p}_3 + n_1 \bmi{p}_2 + n_2 \bmi{p}_1 + n_3 \bmi{p}_0) \Big]\bcdot \hat{\bmi{k}} = 0, \quad \text{at} \quad z = 1
	\end{equation}
\end{subequations}
Using a new variable $\widehat{n}_3=\int_{z}^{1} n_3 \ud z$, the governing system can be written as
\begin{equation}
	\bmi{\mathcal{L}} \bmi{X}_3 = \bmi{F}_3
\end{equation}
\begin{equation}
	\implies
	{\setlength{\arraycolsep}{12pt}
		\renewcommand{\arraystretch}{1.6}
		\begin{bmatrix}
			\bnabla^4 - S_c^{-1} \dfrac{\partial}{\partial \tau} \bnabla^2 & R_c \left(\dfrac{\partial^2}{\partial x^2} + \dfrac{\partial^2}{\partial y^2}\right) \dfrac{\partial}{\partial z} \\
			-\dfrac{\partial n_0}{\partial z} & \dfrac{\partial^2}{\partial \tau \partial z} - \mathcal{F}
	\end{bmatrix}}
	\begin{bmatrix}
		w_3 \\
		\widehat{n}_3
	\end{bmatrix}
	=
	\begin{bmatrix}
		f_{31} \\
		f_{32}
	\end{bmatrix}.
\end{equation}
Linear theory determines the instability threshold and the associated critical mode, but does not describe the finite-amplitude evolution beyond onset. Near criticality, this evolution is governed by nonlinear interactions between the perturbation fields, which enter the resonant problem at third order through
$
\bmi{F}_3=
{\setlength{\arraycolsep}{7pt}
	\renewcommand{\arraystretch}{0.7}
	\begin{bmatrix}
		f_{31} & f_{32}
\end{bmatrix}}^{T}.
$
The forcing $\bmi{F}_3$ contains contributions generated by interactions between the first and second-order fields. At the critical point, the linear operator is singular, since the corresponding homogeneous problem admits a non-trivial solution. The third-order problem therefore possesses a bounded solution only if the resonant forcing satisfies the Fredholm solvability condition
\begin{equation}
	\left\langle \bmi{X}^{\dagger},\bmi{F}_3\right\rangle=0.
\end{equation}
where \(\bmi{X}^{\dagger}\) is the adjoint critical eigenfunction satisfying $\bmi{\mathcal{L}}^{\dagger} \bmi{X}^{\dagger} = 0$. The adjoint operator $\bmi{\mathcal{L}}^{\dagger}$ is defined through $\left\langle\bmi{X}^{\dagger},\bmi{\mathcal{L}} \bmi{X}\right\rangle=\left\langle\bmi{\mathcal{L}}^{\dagger} \bmi{X}^{\dagger},\bmi{X}\right\rangle$. Application of this solvability condition yields the Stuart-Landau amplitude equation
\begin{equation}
	\frac{\ud A}{\ud s} = a A - b |A|^2 A
\end{equation}
which governs the slow evolution of the critical-mode amplitude $A(s)$ in the vicinity of the instability threshold.

\section{Numerical results}\label{4th}
Here, the optical and swimming
parameters are held fixed at $S_c=20$, $V_c=15$, $G_c=1$, $\tau_H=1$ and
$\omega=0.59$, and only the incidence angle is varied, taking
$\theta_i=0^\circ$, $40^\circ$ and $80^\circ$. 

\subsection{Case 1: Stress-free top surface}
At $\theta_i=0^\circ$, the linear stability analysis (first-order problem)
predicts $R_c=314.94$ at $k_c=2.71$. These critical values lie on the
stationary branch of the neutral curve, which implies a stationary
bioconvective solution. The Landau coefficients obtained from the
solvability condition are $a=18.52$ and $b=1.502$. Since $b_r>0$, the
bifurcation is a supercritical pitchfork. When $\theta_i$ is increased to $40^\circ$, the critical Rayleigh number
decreases to $R_c=275.89$ at $k_c=2.90$, and the solution remains
stationary. For this case, the Landau coefficients are $a=23.42$ and
$b=-0.126$. Since $b_r<0$, the bifurcation is a subcritical pitchfork. When $\theta_i$ is further increased to $80^\circ$, the critical Rayleigh
number decreases to $R_c=194.58$, and the critical wavenumber drops to
$k_c=1.16$, so that the pattern wavelength more than doubles. The
solution is again stationary. The neutral curve is very flat near its
minimum, which is reflected in the small value $a=1.785$, and the cubic
coefficient is large and negative, $b=-179.0$. The bifurcation is
therefore strongly subcritical. The dominant contribution to $b$ comes
from the interaction of the critical mode with its second harmonic. This
is because the wavenumber $2k_c$ lies only slightly above the neutral
curve, so that the second harmonic is itself weakly damped. Thus, for a stress-free top surface, increasing the incidence angle
lowers the critical Rayleigh number monotonically and changes the
bifurcation from supercritical to increasingly subcritical.

\subsection{Case 2: Rigid top surface}
At $\theta_i=0^\circ$, the linear stability analysis predicts
$R_c=434.50$ at $k_c=3.41$. These critical values lie on the stationary
branch of the neutral curve, which implies a stationary bioconvective
solution. The Landau coefficients are $a=32.35$ and $b=0.287$. Since
$b_r>0$, the bifurcation is a supercritical pitchfork. When $\theta_i$ is increased to $40^\circ$, the critical Rayleigh number
increases to $R_c=675.95$ at $k_c=3.07$, and the minimum of the neutral
curve now lies on the oscillatory branch, with frequency
$\omega_c=10.28$. The solution therefore becomes overstable
(oscillatory). For this case, the Landau coefficients are
$a=18.83-8.64\,\mathrm{i}$ and $b=0.652+0.673\,\mathrm{i}$. Since
$b_r>0$, the bifurcation is a supercritical Hopf bifurcation, and $b_i$
gives the leading nonlinear correction to the oscillation frequency. When $\theta_i$ is further increased to $80^\circ$, the critical Rayleigh
number increases again, to $R_c=947.76$ at $k_c=3.93$. The critical mode
returns to the stationary branch of the neutral curve, so that the
solution is again stationary. For this case, the Landau coefficients are
$a=31.86$ and $b=-2.946$. Since $b_r<0$, the bifurcation is a subcritical
pitchfork. Thus, for a rigid top surface, increasing the incidence angle raises the
critical Rayleigh number monotonically. The onset changes from
stationary to oscillatory and back to stationary, while the bifurcation
changes from a supercritical pitchfork to a supercritical Hopf
bifurcation and finally to a subcritical pitchfork.

\begin{table}[!htbp]
	\centering
	\caption{Critical states and Stuart--Landau coefficients
 for the case of $S_c=20$, $V_c=15$, $G_c=1$, $\tau_H=1$, $\omega=0.59$. Here SF: stress-free upper surface;
		R: rigid upper surface; SC: supercritical; sub: subcritical.}
	\label{tab:fourcases}
	{\small
		\begin{tabular}{lccccccc}
			\toprule
			Case & $\theta_i$ & $R_c$ & $k_c$ & $\omega_c$ &
			$a$ & $b$ & bifurcation\\
			\midrule
			SF  & $0^\circ$  & 314.94 & 2.706 & 0     & 18.52 & $1.502$  & SC pitchfork \\
			SF & $40^\circ$ & 275.89 & 2.900 & 0     & 23.42 & $-0.126$ & sub pitchfork \\
			SF & $80^\circ$ & 194.58 & 1.163 & 0     & 1.785 & $-179.0$ & sub pitchfork \\
			\midrule
			R   & $0^\circ$  & 434.50 & 3.408 & 0     & 32.35 & $0.287$  & SC pitchfork  \\
			R  & $40^\circ$ & 675.95 & 3.075 & 10.28 & $18.83-8.64\,\mathrm{i}$ & $0.652+0.673\,\mathrm{i}$ & SC Hopf \\
			R  & $80^\circ$ & 947.76 & 3.932 & 0     & 31.86 & $-2.946$ & sub pitchfork \\
			\bottomrule
	\end{tabular}}
\end{table}
\section{Conclusions}\label{5th}

The weakly nonlinear evolution of phototaxis-induced bioconvection has been studied for an absorbing and isotropically scattering algal suspension illuminated by oblique collimated irradiation, with normal incidence recovered as the limit $\theta_i=0^\circ$. The lower boundary is rigid, and both stress-free and rigid upper surfaces are considered. A multiple-scales expansion to third order, together with the Fredholm solvability condition, leads to a Stuart--Landau equation whose coefficients determine the character of the bifurcation near onset. The main conclusions are as follows.

\begin{enumerate}
	\item Oblique irradiation acts on the onset threshold in opposite directions for the two configurations. As the incidence angle increases, the critical Rayleigh number decreases monotonically beneath a stress-free surface, and increases monotonically beneath a rigid surface.
	
	\item For a stress-free upper surface, the most unstable mode is stationary at all incidence angles considered. At normal incidence the bifurcation is a supercritical pitchfork, and the bioconvective rolls saturate smoothly. With increasing obliquity the real part of the cubic Landau coefficient changes sign and the bifurcation becomes subcritical. At strongly oblique incidence the subcriticality becomes pronounced and the preferred wavelength of the pattern increases considerably. 
	
	\item For a rigid upper surface, the critical Rayleigh number is higher than for a stress-free surface at every incidence angle. At normal incidence the bifurcation is a supercritical pitchfork. At moderate obliquity the critical mode moves to the oscillatory branch of the neutral curve, and the resulting overstable solution arises through a supercritical Hopf bifurcation. The imaginary part of the Landau coefficient gives the nonlinear correction to the
	oscillation frequency. At strong obliquity the critical mode returns to the stationary branch, and the pitchfork bifurcation becomes subcritical.
	
	\item In both configurations, strongly oblique irradiation leads to
	a subcritical bifurcation. Finite-amplitude bioconvection may then
	exist below the linear threshold.
\end{enumerate}


\begin{appendices}
	\numberwithin{equation}{section}
	
	\section{Perturbation expansion of the phototactic fields}
	\label{app:light}
	
	The phototactic response is expanded as
	\begin{align}
		\nonumber\bmi{p} &= \bmi{p}_0 + \epsilon \bmi{p}_1 + \epsilon^2 \bmi{p}_2 + \epsilon^3 \bmi{p}_3 + \cdots \\
		& = -T(G_0 + \epsilon G_1 + \epsilon^2 G_2 + \epsilon^3 G_3 + \cdots) \frac{\bmi{q}_0 + \epsilon \bmi{q}_1 + \epsilon^2 \bmi{q}_2 + \epsilon^3 \bmi{q}_3 + \cdots}{|\bmi{q}_0 + \epsilon \bmi{q}_1 + \epsilon^2 \bmi{q}_2 + \epsilon^3 \bmi{q}_3 + \cdots|}
	\end{align} 
	
	After solving, we can find the value of 
	$\bmi{p}_0, \bmi{p}_1, \bmi{p}_2, \bmi{p}_3 $
	
\end{appendices}

\section*{Data availability}
The data presented in this article support the findings of this study.
\nocite{*}


\bibliographystyle{plainnat}
\bibliography{skrajput}

\end{document}